\documentclass[11pt]{amsart}

\usepackage{amsmath, amssymb, amscd, graphicx}

\usepackage[all]{xy}

\newtheorem{thm}{Theorem}[section]
\newtheorem{conj}[thm]{Conjecture}

\newtheorem{lem}[thm]{Lemma}
\newtheorem{prop}[thm]{Proposition}

\newtheorem{defin}[thm]{Definition}

\def\q{{\mathcal Q}}

\def\R{{\mathbb R}}

\def\Z{{\mathbb Z}}

\def\del{{\partial}}

\def\mod{{\textup{mod} \;}}

\begin{document}

\title[H-thin, non-QA links]%
{Homologically thin, non-quasi-alternating links}

\begin{abstract}

We exhibit the first examples of links which are homologically thin but not quasi-alternating.  To show that they are not quasi-alternating, we argue that none of their branched double-covers bounds a negative definite 4-manifold with non-torsion $H_1$.  Using this method, we also complete the determination of the quasi-alternating pretzel links.

\end{abstract}

\author[Joshua Greene]{Joshua Greene}

\address{Department of Mathematics, Princeton University\\ Princeton, NJ 08542}

\thanks{Partially supported by an NSF Graduate Fellowship.}

\email{jegreene@math.princeton.edu}

\maketitle

\section{Introduction.}

Quasi-alternating links were defined by Ozsv\'ath and Szab\'o \cite[Definition 3.9]{os:doublecover}. They are a natural generalization of the class of alternating links.

\begin{defin}\label{def: qa} The set $\q$ of {\em quasi-alternating (QA) links} is the smallest set of links such that

\begin{itemize}

\item the unknot $U$ belongs to $\q$, and

\item if $L$ is a link with a projection containing a crossing for which the two resolutions $L_0$ and $L_1$ belong to $\q$, and $\det(L) = \det(L_0) + \det(L_1)$, then $L$ belongs to $\q$.

\end{itemize}

\end{defin}

By this definition, non-split alternating links belong to $\q$.  Many familiar properties of alternating links hold for a QA link $L$:

\begin{enumerate}

\item the branched double-cover $\Sigma(L)$ is an L-space \cite[Proposition 3.3]{os:doublecover};

\item the space $\Sigma(L)$ bounds a negative definite 4-manifold $W$ with $H_1(W) = 0$ \cite[Proof of Lemma 3.6]{os:doublecover};

\item the $\Z/2\Z$ knot Floer homology group $\widehat{HFK}(L;\Z/2\Z)$ is thin \cite[Theorem 2]{mo:qa};

\item the reduced ordinary Khovanov homology group $\overline{Kh}(L)$ is thin \cite[Theorem 1]{mo:qa}; and

\item the reduced odd Khovanov homology group $\overline{Kh}'(L)$ is thin \cite[Remark after Proposition 5.2]{ors:oddkh}.

\end{enumerate}

It is an interesting open problem to characterize those links that are homologically thin with respect to any of the above knot homology theories.  For some time, it remained a possibility that a link was $\widehat{HFK}$- or $\overline{Kh}$-thin if and only if it was QA.  This possibility was recently refuted by Shumakovitch, who used his excellent program KhoHo \cite{s:khoho} to show that the pretzel knots $9_{46} = P(3,-3,3)$ and $10_{140} = P(3,-3,4)$ have torsion in their odd Khovanov homology groups, although they are both $\widehat{HFK}$- and $\overline{Kh}$-thin.  Thus, neither of these knots is {\em odd-thin}, so neither is QA.\footnote{In fact, in all known examples, an odd-thin link is $\overline{Kh}$-thin, and a link is $\overline{Kh}$-thin iff it is $\widehat{HFK}$-thin.  A conjecture of Rasmussen would imply that a $\overline{Kh}$-thin link is necessarily $\widehat{HFK}$-thin \cite[Section 5]{r:knothoms}.}

\begin{defin}\label{def: thin}

A link $L$ is {\em homologically thin} (without qualification) if it is simultaneously thin with respect to $\widehat{HFK}, \overline{Kh},$ and $\overline{Kh}'$.

\end{defin}

It has remained a challenge to exhibit a link that is homologically thin and not QA.  The purpose of this note is to describe such examples, and moreover to exploit the condition (2) as an obstruction to QA-ness.

\begin{thm}\label{thm: non-qa}

There exist homologically thin, non-QA links.

\end{thm}

At the heart of our method is Donaldson's celebrated ``Theorem A'', which asserts that the intersection pairing of a smooth, closed, negative definite 4-manifold is diagonalizable \cite{d:thma}.  Coupled with calculations by several researchers \cite{bald:3braids,ck:qa,js:calc,m:triangle,s:khoho,widmer:qa}, we identify $11n50$ as the only knot with up to 11 crossings which is neither QA nor odd-thick.  Furthermore, combined with work of Champanerkar-Kofman \cite[Theorem 3.2]{ck:qa}, we complete the determination of the QA pretzel links.  For a clear, concise account of the relevant notation and facts concerning Montesinos links here and in what follows, see \cite[Section 3.2]{os:qhs}.

\begin{figure}
\centering
\includegraphics[width=1.5in]{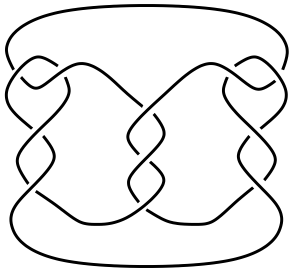}
\caption{The knot $11n50$.}  \label{f: 11n50}
\end{figure}

\begin{thm}\label{thm: pretzel}

The pretzel link $P(-e; p_1,\dots,p_n, -q_1,\dots ,-q_m) = M(-e;(p_1,1),\dots,(p_n,1),$ $(q_1,-1),\cdots,(q_m,-1))$, with $e, n, m \geq 0$, all $p_i \geq 2$, and all $q_j \geq 3$, is QA iff

\begin{enumerate}

\item $e > m - 1$;

\item $e = m - 1 > 0$;

\item $e = 0, n = 1,$ and $p_1 > \min \{ q_1, \dots, q_m \}$ or $m \leq 1$; or

\item $e = 0, m = 1,$ and $q_1 > \min \{ p_1, \dots, p_n \}$ or $n \leq 1$.

\end{enumerate} The same is true on permuting the parameters $p_i$ and $q_j$.

\end{thm}

Any pretzel link can be put in the above form after mirroring \cite[Section 2.3]{kawauchi}, which clearly preserves the QA property.  Section \ref{s: main} contains the proofs of Theorems \ref{thm: non-qa} and \ref{thm: pretzel}.  Section \ref{s: discussion} contains some related examples, as well as some discussion surrounding Conjecture \ref{conj: qa det}, which asserts that there are finitely many QA links of bounded determinant.



\section*{Acknowledgments.}

Thanks to John Baldwin, Abhijit Champanerkar, Michael Eisermann, Slavik Jablan, and Liam Watson for helpful correspondence.


\section{Proofs of the main results.}\label{s: main}

To prove the theorems, we rely on the following lemma.  We use (co)homology groups with integer coefficients throughout.

\begin{lem}\label{l:key}

Suppose that $X$ and $W$ are a pair of 4-manifolds, $\del X = - \del W = Y$ is a rational homology sphere, and $H_1(W)$ is torsion-free. Express the map $H_2(X)/\textup{Tors} \to H_2(X \cup W)/\textup{Tors}$ with respect to a pair of bases by the matrix $A$.  This map is an inclusion, and if some $k$ rows of $A$ contain all the non-zero entries of some $k$ of its columns, then the induced $k \times k$ minor has determinant $\pm 1$.

\end{lem}


\begin{proof}

The stated assumption ensures that the restriction map $H^2(X \cup W) \to H^2(X)$ surjects.  Just to be sure, consider the long exact sequences in cohomology for the pairs $(X \cup W,X)$ and $(W,Y)$, and the natural map between them.  The relevant portion reads \[
\xymatrix{ H^2(X \cup W) \ar[r] & H^2(X) \ar[r] \ar[d] & H^3(X \cup W, X) \ar[d] \\  & H^2(Y) \ar[r] & H^3(W,Y). }
\] The second vertical arrow is an isomorphism by excision, and Poincar\'e-Lefschetz duality identifies this group with $H_1(W)$, which is torsion-free.  Since $H^2(Y)$ is torsion, the bottom horizontal map is $0$.  It follows that the map $H^2(X) \to H^3(X \cup W, X)$ is 0, so the map $H^2(X \cup W) \to H^2(X)$ surjects as claimed.



Consequently, the map $H^2(X \cup W)/ \text{Tors} \to H^2(X) / \text{Tors}$ surjects as well.  On the other hand, this map of groups is dual to the map $H_2(X)/\text{Tors} \to H_2(X \cup W)/\text{Tors}$, so is given with respect to the pair of dual bases by the matrix $A^T$.  Suppose that some $k$ rows of $A$ contain all the non-zero entries of some $k$ of its columns, and let $B$ denote the corresponding $k \times k$ minor.  By permuting the basis elements if necessary, we may assume that $B$ is the top-left $k \times k$ minor, possibly changing its determinant by a sign:  \[A = \left( \begin{array}{cc} B & C \\ 0 & D \end{array} \right).\]  Since the dual map $A^T$ surjects, the map $B^T$ must as well, hence $\det(B) = \pm 1$, as claimed.  The fact that the map $A$ injects follows, for example, from the Mayer-Vietoris sequence for the natural decomposition of $X \cup W$, noting that $H_2(Y)$ vanishes.

\end{proof}

\begin{proof}[Proof of Theorem \ref{thm: non-qa}]

We establish the result by showing that $K =$ $11n50$ is homologically thin but not QA.  Additional examples appear in Subsection \ref{ss: examples}.  The knot Floer homology group $\widehat{HFK}(K;\Z/2\Z)$ was calculated by Baldwin-Gillam \cite{bg:hfk}, and the Khovanov homology groups $\overline{Kh}(K)$ and $\overline{Kh}'(K)$ by Shumakovitch \cite{s:khoho}; their results show that $K$ is homologically thin.

Figure \ref{f: 11n50} exhibits the knot $11n50$ as the Montesinos knot $M(0;(5,2),(3,1),(5,-2))$, which is equivalent to $M(1;(5,2),(3,1),(5,3))$. As such, its branched double-cover $Y$ is the boundary of the plumbing $X$ on the graph $\Gamma$ shown in Figure \ref{f: plumbing}.  We label the vertices of $\Gamma$ from left to right and top to bottom by $v_1,\dots,v_7$.  The spheres corresponding to these vertices give rise to a basis for $H_2(X)$, with respect to which the intersection pairing on $X$ is given by the weighted adjacency matrix $A_\Gamma$, whose $(i,i)$-entry records the weight of $v_i$, and whose $(i,j)$-entry for $i \ne j$ is $1$ or $0$ according as $v_i$ and $v_j$ are adjacent or not.  The space $X$ is negative definite and $H_1(X) = 0$.

\begin{figure}
\centering
\includegraphics[width=2in]{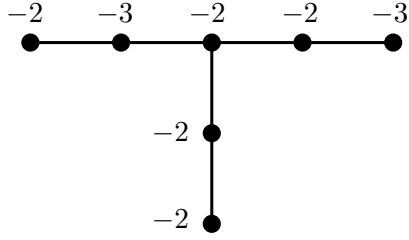}
\put(-150,80){$-2$}
\put(-116,80){$-3$}
\put(-81,80){$-2$}
\put(-46,80){$-2$}
\put(-12,80){$-3$}
\put(-95,35){$-2$}
\put(-95,2){$-2$}
\caption{The graph $\Gamma$.}  \label{f: plumbing}
\end{figure}

If $11n50$ were QA, then its mirror $\overline{11n50}$ would be as well, whose branched double-cover is $-Y$.  We proceed by way of fact (2) to derive a contradiction.  According to it, there must exist a negative definite 4-manifold $W$ with $\del W = -Y$ and $H_1(W) = 0$.  Consider the 4-manifold $X \cup W$.  It is a closed, smooth, negative definite 4-manifold, so by Donaldson's theorem, its intersection pairing is diagonalizable. That is, there exists an isomorphism $(H_2(X \cup W)/\text{Tors},Q_{X \cup W}) \cong - \Z^n$, where $n = b_2(X \cup W)$ and $-\Z^n = \langle E_1,\dots,E_n \rangle$ denotes the space $\Z^n$ equipped with minus its standard positive definite inner product.

Therefore, $(H_2(X)/\text{Tors},Q_X)$ embeds into $-\Z^n$, for some $n$.  Let $x_i$ denote the image of the class in $H_2(X)/\text{Tors}$ under the inclusion into $-\Z^n$ which corresponds to the vertex $v_i$.  No pair of vertices of weight -2 has the same set of neighbors, so the corresponding vectors in $-\Z^n$ have distinct reductions $(\mod 2)$.  This observation helps us to deduce that, by applying an automorphism to $-\Z^n$, we have $x_7 = E_1-E_2, x_6 = E_2-E_3, x_3 = E_3-E_4, x_4 = E_4-E_5$, and $x_1 = E_6-E_7$. Now, swapping $E_6$ and $-E_7$ if need be, we obtain $x_2 = E_4+E_5-E_6$, and then $x_5 = E_5+E_6+E_7$.  Thus, with respect to the chosen bases for $H_2(X)/\text{Tors}$ and $H_2(X \cup W)/\text{Tors}$, the inclusion map is given by a matrix $A$ whose seven columns are supported in its first seven rows.  Let $B$ denote the induced $7 \times 7$ minor.  Then $-B^T B = -A^T  = A_\Gamma$, and this is a presentation matrix for $H^2(Y) \cong \Z / 25 \Z$.
Hence $|\det(B)| = 5 \ne 1$, in contradiction to Lemma \ref{l:key}.

It follows that $-Y$ does not bound a negative definite 4-manifold with torsion-free, let alone vanishing, $H_1$, and so the knot $11n50$ is not QA.

\end{proof}

We remark that there does exist a negative definite 4-manifold $W$ with boundary $-Y$ for which $H_1(W)$ contains torsion.  The knot $11n50$ is a slice knot\footnote{Michael Eisermann points out that this is direct from the presentation of $11n50$ as a {\em symmetric union} in Figure \ref{f: 11n50} (cf. \cite[Theorem 5]{lamm:ribbons}).}, so the double-cover of $D^4$, branched along a slice disk for the mirror $\overline{11n50}$, is a rational homology ball with boundary $-Y$, which we may blow up to make undeniably negative definite.  However, its first homology group contains a subgroup isomorphic to a quotient of $H_1(Y)$ of square-root order, which is $\Z/ 5\Z$ in this case.

\begin{prop}\label{p: pretzel}

For $n \geq 2$, and $p_1,\dots,p_n \geq 2$, and $q \geq 1$, the pretzel link $P(p_1,\dots,p_n,-q)$ is QA iff $q > \min \{p_1,\dots,p_n\}$.

\end{prop}

\begin{proof}

Let $L$ denote the pretzel link $P(p_1,\dots,p_n,-q) = M(1;(p_1,1),\dots,(p_n,1),(q,q-1))$.  The space $\Sigma(L)$  is the boundary of the plumbing $X$ on a star-shaped planar graph $\Gamma$.  The graph $\Gamma$ has $n+1$ rays of lengths $p_1-1,\dots,p_n-1$, and $1$ emanating from the star vertex in cyclic order; by {\em length} we mean the number of edges.  The star vertex receives weight $-n$, the vertex on the distinguished ray of length $1$ receives weight $-q$, and every other vertex receives weight $-2$.  As before, the intersection pairing on $X$ is given in the natural spherical basis by the weighted adjacency matrix $A_\Gamma$.

The space $X$ is negative definite if and only if  $p_1^{-1}+\cdots+p_n^{-1}-q^{-1} > 0$ \cite[Theorem 5.2]{nr:seifert}.  If it is not -- implying that $q < \min \{ p_1,\dots, p_n \}$ -- then we claim that $L$ is not QA.  Consider the space $-Y = \Sigma(\overline{L})$.  It is the boundary of plumbing on a star-shaped graph $\Gamma'$ with $n$ rays of length $1$ and one of length $q-1$ emanating from the star vertex.  The vertices on rays of length $1$ receive weights $-p_1,\dots,-p_n$ in cyclic order, the vertices on the ray of length $q-1$ receive weight $-2$, and the star vertex receives weight $-1$.  Under the assumption that $\Gamma$ is {\em not} negative definite, the graph $\Gamma'$ {\em is} by another application of \cite[Theorem 5.2]{nr:seifert}.


Now we appeal to some facts for which a detailed account would extend too far beyond the scope of this note.  Since the star vertex has degree $\geq 3$ and weight $-1$, an application of Laufer's algorithm terminates at the $0^{th}$ iteration, and shows at once that the space $-Y = Y_{\Gamma'}$ is not the link of a rational surface singularity \cite[Section 4]{laufer:rational}.  The invariant $HF^+(Y) = HF^+(-Y_{\Gamma'})$ is identified with a particular $\Z[U]$-module $\mathbb{H}^+(\Gamma')$, as detailed in \cite[Section 2]{os:plumbed}.  On the other hand, N\'emethi has proven that if $HF^+(-Y_{\Gamma'}) \cong \mathbb{H}^+(\Gamma')$ and $Y_{\Gamma'}$ is an L-space, then $Y_{\Gamma'}$ is the link of a rational surface singularity \cite[Proposition 4.1.2]{n:latticecohomology}.  It follows that $Y$ is not an L-space, so $L$ is not QA in this case.

Therefore, we may assume henceforth that $X$ is negative definite.  Now suppose that $L$ were QA, so that $-\Sigma(L) = \del W$, where $W$ is a negative definite 4-manifold with $H_1(W) = 0$.  Since $\det(L) > 0$, it follows that $-\Sigma(L)$ is a rational homology sphere. We proceed as in the proof of Theorem \ref{thm: non-qa} and analyze how $(H_2(X),Q_X)$ can embed into the lattice $-\Z^n = \langle E_1,\dots,E_n \rangle$.  To every vertex of $\Gamma$ corresponds a vector in $-\Z^n$.  If two distinct vertices of weight $-2$ gave rise to vectors with the same reduction $(\mod 2)$, then a change of basis of $-\Z^n$ puts these vectors in the form $E_1+E_2$ and $E_1-E_2$.  These in turn correspond to a pair of columns of the matrix $A$ representing the map $H_2(X)/\text{Tors} \to H_2(X \cup W)/\text{Tors}$ supported in the first two columns. The induced $2 \times 2$ minor has determinant $\pm 2$, in contradiction to Lemma \ref{l:key}.


It follows that, moving away from the star vertex, the vectors corresponding to the vertices along a ray of length $p_i-1$ can be put in the form $E^i_1-E^i_2,\dots,E^i_{p_i-1}-E^i_{p_i}$, where all the basis vectors $E^i_j$, $1 \leq i \leq n$, $1 \leq j \leq p_i$, are distinct.  It then follows that the star vertex corresponds to the vector $-E^1_1 - \cdots - E^n_1$.  Consider the vector $x$ corresponding to the vertex of weight $-q$.  Its inner product with $-E^1_1 - \cdots - E^n_1$ is non-zero, so in its expansion with respect to the chosen basis of $-\Z^n$, it has some term of the form $a \cdot E^i_1$ with $a \ne 0$.  Since $x$ is orthogonal to all those vectors corresponding to the $i^{th}$ ray of $\Gamma$, its expansion takes the form $a (E^i_1 + \cdots + E^i_{p_i}) \; + $ (additional terms).  It follows that $q \geq |a| \cdot p_i \geq \min \{ p_1,\dots,p_n \}$.  Now suppose by way of contradiction that equality held throughout.  Then in fact $x = E^i_1 + \cdots + E^i_{p_i}$.  Consider the rows of the matrix $A$ corresponding to the vectors $E^i_1-E^i_2, \dots, E^i_{p_i-1}-E^i_{p_i}$, and $E^i_1 + \cdots + E^i_{p_i}$.  These are $p_i$ columns supported in $p_i$ rows, whose induced $p_i \times p_i$ minor has determinant $\pm p_i \ne \pm 1$, in contradiction to Lemma \ref{l:key}.  Consequently, $q > \min \{ p_1,\dots,p_n \}$, as desired.

The converse statement in \cite[Theorem 3.2(1)]{ck:qa} completes the proof of the Proposition.

\end{proof}

\begin{proof}[Proof of Theorem \ref{thm: pretzel}]

Let $L$ denote the pretzel link appearing in the statement of the Theorem.  As a Montesinos link, it is notated by $M(e;(p_1,1),\dots,(p_n,1),(q_1,-1),\dots,(q_m,-1))$.

If $e < m-1$, then $L$ has the equivalent description $M(0;(p_1,1),\dots,(p_n,1),(q_1,q_1-1),\dots,(q_e,q_e-1),(q_{e+1},-1),\dots,(q_m,-1))$.  If $e + n \geq 2$, then the diagram resulting from this description is adequate and non-alternating.  It follows that $L$ is $\overline{Kh}$-thick \cite[Proposition 5.1]{khovanov:patterns} and hence not QA in this case.  If $e + n \leq 1$, then either $e = 1, n = 0$, or $e = 0, n = 1$.  In the first case, $L$ is QA if $m \leq 1$ (falling under Case (1) of the Theorem) and non-QA if $m \geq 2$ (applying Proposition \ref{p: pretzel} to $\overline{L}$).  In the second case, Proposition \ref{p: pretzel} applies once again to $\overline{L}$ to determine when $L$ is QA (Case (3)).  This establishes the Theorem in case $e < m - 1$.

If $e > m - 1$, then $L$ has the equivalent description $M(e-m;(p_1,1),\dots,(p_n,1),(q_1,q_1-1),\dots,(q_m,q_m-1))$.  Its associated diagram is connected and alternating, so $L$ is QA in this case (Case (1)).  Also, if $e = m - 1 = 0$, then the Theorem follows by a combination of Proposition \ref{p: pretzel} and \cite[Theorem 3.2(2)]{ck:qa} (Case (4)).

It stands to consider the case that $e = m - 1 > 0$ (Case (2)).  We prove that $L$ is QA by induction on $e+q_1+\cdots+q_m$.  Consider a crossing appearing in the tassle with $-q_m$ half-twists.  The resolution $L_0$ is the link $P(e;p_1,\dots,p_n,-q_1,\dots,-q_{m-1})$, while the resolution $L_1$ is the link $P(e;p_1,\dots,p_n,-q_1,\dots,-(q_m-1))$.  We calculate \[
\begin{array}{rcl}
\det(L_0) &=& p_1 \cdots p_n q_1 \cdots q_{m-1} (e + p_1^{-1} + \cdots + p_n^{-1} - q_1^{-1} - \cdots - q_{m-1}^{-1}), \cr \\
\det(L_1) &=& p_1 \cdots p_n q_1 \cdots q_{m-1} (q_m - 1) (e + p_1^{-1} + \cdots + p_n^{-1} - q_1^{-1} - \cdots - q_{m-1}^{-1} - (q_m-1)^{-1}), \cr \\
\det(L) &=& p_1 \cdots p_n q_1 \cdots q_m (e + p_1^{-1} + \cdots + p_n^{-1} - q_1^{-1} - \cdots - q_m^{-1}).
\end{array} \] Note in particular that the expression for each determinant is positive, since $e = m - 1$ and there are at most $m$ negative terms in each sum, with each term $\geq -1/2$.  The equality $\det(L) = \det(L_0) + \det(L_1)$ is immediate.  Now, $L_0$ has the same value $e$ and one fewer negative term, so as in the case $e > m-1$ treated above, this link has a connected, alternating diagram, hence is QA.  If $q_m > 3$, then the link $L_1$ is QA by induction.  Otherwise, $q_m = 3$, and so $L_1 = P(e;p_1,\dots,p_n,-q_1,\dots,-q_{m-1},-2) = P(e-1;p_1,\dots,p_n,-q_1,\dots,-q_{m-1},2) = P(e-1;2,p_1,\dots,p_n,-q_1,\dots,-q_{m-1})$.  If $e-1 = m-2 > 0$, then $L_1$ is QA by induction, while if $e-1 = m-2 = 0$, then $L_1$ is QA by Proposition \ref{p: pretzel}.  Thus, $L_1$ is QA regardless, and it follows that $L$ is QA as well.  This completes the induction step.

The preceding argument carries over {\em mutatis mutandis} to the case of a pretzel link which differs from $L$ by a permutation of the parameters $p_i$ and $q_j$.  This completes the proof of the Theorem.

\end{proof}


\section{Discussion.}\label{s: discussion}

\subsection{Further obstructions.}

The main principle at work in this note is the fact that for a QA link $L$, there is naturally associated to it a smooth, negative definite 4-manifold $X_L$ with vanishing $H_1$ and boundary $\Sigma(L)$.  It is therefore of interest to have on hand obstructions to a 3-manifold bounding a negative definite 4-manifold with torsion-free or vanishing $H_1$, and to examine more closely the topology of $X_L$ in the hopes of developing finer obstructions to QA-ness.

In the first direction, Ozsv\'ath-Szab\'o \cite[Section 9.2]{os:absgr} have developed an obstruction which makes use of the {\em correction terms} in Heegaard Floer homology, and which was subsequently sharpened by Owens-Strle \cite[Theorem 2]{os:lattices}.    For the case of $-Y = \Sigma(\overline{11n50})$, the Owens-Strle obstruction does not rule out the possibility that this space bounds a negative definite 4-manifold with torsion-free $H_1$.  Indeed, using the plumbing graph $\Gamma$, we can calculate the correction terms of $-Y$ according to \cite[Corollary 1.5]{os:plumbed}.  The largest correction term has the value $8/25$, which passes their obstruction since it is $> 1/4$.  Therefore, the argument given in Theorem \ref{thm: non-qa} provides information where this obstruction does not.

In the second direction, Ozsv\'ath-Szab\'o have shown that $X_L$ is a {\em sharp} 4-manifold when $L$ is an alternating link (\cite[Section 2.8]{os:unknotting}, \cite[Theorem 3.4]{os:doublecover}), and an early arXiv version of the paper \cite{os:doublecover} suggested that the same is true for an arbitrary QA link $L$ ({\tt math.GT/0309170}, after Proposition 3.3).  However, this is not the case.  Indeed, $L = \overline{8_{20}} = P(3,-2,-2)$ does not bound {\em any} sharp 4-manifold \cite[Proposition 7.3]{g:3braids}.  This negative result begs for an efficient means of calculating the correction terms of $\Sigma(L)$ for a QA link $L$ in general.  Is it still possible to utilize $X_L$ in some way towards this end?  What further information can we glean from $X_L$ to develop an obstruction to QA-ness?

\subsection{Further examples.}\label{ss: examples}

The connect-sum of $11n50$ with any QA link $L$ will result in a homologically thin link which is not QA.  The fact that it is homologically thin follows from the behavior of the relevant knot homology groups under the connect-sum operation.  The fact that it is not QA follows the proof of Theorem \ref{thm: non-qa}, noting that $\Sigma(11n50 \# L)$ is the boundary of $X_{11n50} \# X_L$, whose intersection pairing decomposes as a direct sum.

In the interest of giving examples of {\em prime} links that are homologically thin but not QA, consider the family of examples given by the prime Montesinos links $L(m,n) = M(0;(m^2+1,m),(n,1),(m^2+1,-m))$ for positive integers $m,n \geq 2$.  Thus, $11n50$ is the knot $L(2,3)$.  The proof of Theorem \ref{thm: non-qa} easily generalizes to show that for all $n > m$, the link $L(m,n)$ is not QA.  Jablan-Sazdanovi\'c \cite{js:calc} suggested the family $L(2,n)$ with $n \geq 3$ as an infinite family of non-QA, homologically thin links.  Indeed, an application of the skein sequence in knot Floer homology shows that $L(2,n)$ is $\widehat{HFK}$-thin for all $n \geq 6$, and a calculation in KhoHo shows that both the relevant Khovanov homology groups are thin for $n = 3,4,5$.  However, $\overline{Kh}' (L(2,6))$ possesses 5-torsion, and this persists for all $L(2,n), n \geq 6$, by an application of the long exact sequence in Khovanov homology.  Indeed, small calculations suggest that $\overline{Kh}'(L(m,n))$ will possess torsion of order $m^2+1$ for all $n \gg m$.  Furthermore, the pretzel link $P(q,-q,q)$ has $q$-torsion for all $q \leq 6$ \cite{s:khoho,s:oddho}.  It would be very interesting to explain these torsion phenomena.  However, it is reasonable to conjecture that the links $L(m,m+1)$ are homologically thin for all $m \geq 2$, and thereby provide an infinite source of examples of prime, homologically thin, non-QA links.  We discuss another potential family at the end of this note.

The family of links $L(m,n)$, together with work of Champanerkar-Kofman \cite{ck:qa} and Widmer \cite{widmer:qa}, indicate some progress in extending Theorem \ref{thm: pretzel} to the more general case of Montesinos links.  We hope to address this question more fully in future work.


Lastly, extensive calculations by Jablan-Sazdanovi\'c \cite{js:comm,js:calc} (including corrections to some of the ones that originally appeared in \cite{js:calc}) indicate that amongst multi-component {\em links} with up to 11 crossings, all except L11n77 and L11n90 are either odd-thick or non-QA.  The method described here can be used to show that L11n90 is not QA, although it is odd-thin.  We were unable to conclude anything further about L11n90.  Therefore, it may require additional ideas to prove that it is non-QA, if indeed this is the case.


\subsection{A conjecture.}  We close with a conjecture.

\begin{conj}\label{conj: qa det}

There exist only finitely many QA links with a given determinant.

\end{conj}  In support of Conjecture \ref{conj: qa det}, note that there are finitely many alternating links with a given determinant.  Moreover, we have the following result for very small determinants.

\begin{prop}\label{p: det < 4}

If $L$ is a QA link with determinant $\leq 3$, then $L$ is alternating.

\end{prop}

\begin{proof}

Suppose that $L$ is a QA link.  If $\det(L) = 1$, then the assertion is trivial.

Next, suppose that $\det(L) =2$.  Let $c$ denote a QA crossing in a diagram of $L$, and $L_0$ and $L_1$ the two resolutions of $L$ at $c$.  Of course, both $L_0$ and $L_1$ are unknots.  Let $\gamma$ denote a small unknotted arc connecting the two strands nearby the resolution in $L_0$, and let $K$ denote its preimage in $\Sigma(L_0) = S^3$.  Then $\Sigma(L_1) = S^3$ is a non-trivial surgery on $K$; by the Dehn surgery characterization of the unknot (\cite[Theorem 1.1]{kmos}, or \cite[Theorem 2]{gl:knotcomplements} in this special case), it follows that $K$ is the unknot.  The space $\Sigma(L)$ is an integer surgery on $K$ as well, and since $\det(L)=2$ it follows that $\Sigma(L) \cong \R P^3$.  A result of Hodgson-Rubinstein \cite[Corollary 4.12]{hr:2bridge} characterizes $2$-bridge links as those links whose branched double-covers are lens spaces.  It follows that $L$ is the Hopf link.

Lastly, suppose that $L$ is QA and $\det(L)=3$, and proceed as above.  In this case, $L_0$ is the unknot, while $L_1$ is the Hopf link (or vice versa).  Now $\Sigma(L_1) \cong \R P^3$ is a surgery on $K \subset S^3$, and the Dehn surgery characterization of the unknot once again shows that $K$ is the unknot.  Hence $\Sigma(L)$ is the lens space $\pm L(3,1)$, and citing \cite[Corollary 4.12]{hr:2bridge} again shows that $L$ is a trefoil knot.

John Baldwin gives an alternative argument for the case of a determinant $3$ QA link.  Such a link is $\widehat{HFK}$-thin.  Now using the facts that $\widehat{HFK}(L)$ is the $E_1$ term in a spectral sequence converging to $\widehat{HF}(S^3) \cong \Z_{(0)}$, and that $\widehat{HFK}_d(L,i) \cong \widehat{HFK}_{d-2i}(L,-i)$ \cite{os:hfk}, it follows that $L$ has the knot Floer homology of a trefoil knot.  Since a trefoil is uniquely determined by its knot Floer homology \cite{ghiggini}, the result follows.

\end{proof}

We remark that if Conjecture \ref{conj: qa det} were false, and there were infinitely many QA links of some fixed determinant, then amongst their branched double-covers we would obtain an infinite family of L-spaces with the same graded Heegaard Floer homology groups.  No such family is known as of this writing.  The details of this argument will appear in a separate paper.

In contrast to Conjecture \ref{conj: qa det}, Liam Watson points out that there exist infinitely many homologically thin knots with the same homological invariants as the knot $11n50$ (in particular, determinant 25), which we now describe.  The knot $11n50$ occurs as $K(0,3)$ in Kanenobu's two-parameter family of knots $K(p,q)$ whose HOMFLY polynomial depends only on $p+q$ \cite{kanenobu:samepoly}.  Watson showed that the unreduced ordinary Khovanov homology of $K(p,q)$ depends only on the value $p+q$  as well \cite[esp. Lemma 3.1 and Section 7.4]{watson:samekh}, and the same is true in the context of unreduced odd Khovanov homology by a similar argument.  Since the unreduced groups are thin for $K(0,3)$, they are thin and equal for all $K(n,3-n)$, hence the reduced versions are thin and equal as well.  Furthermore, an application of the long exact sequence in knot Floer homology implies that $\widehat{HFK}(K(n,3-n))$ is independent of $n$.  Finally, Kanenobu's work implies that the knots $K(n,3-n)$ are distinguished by their Alexander modules.  Therefore, the knots $K(n,3-n)$ provide the desired source of examples.  Amongst them, the knot $K(1,2) = K(2,1) = 11n132$ is QA.  However, we speculate that this is the only knot in this family that is QA, and this is the subject of work in progress.




\bibliographystyle{plain}
\bibliography{References}

\end{document}